\documentclass[12pt,reqno]{amsart}
\usepackage{amssymb,amsfonts,amsmath,amsthm,charter}

\usepackage{orcidlink}

\usepackage[a4paper,margin=2.5cm,top=2.5cm,bottom=2.5cm,centering,headheight=3ex,headsep=4ex,vcentering]{geometry}

\definecolor{ao(english)}{rgb}{0.13, 0.55, 0.13}

\hypersetup{
  colorlinks=true,
  linkcolor=ao(english),
  citecolor=ao(english),
  urlcolor=ao(english)
}

\theoremstyle{plain}

\numberwithin{equation}{section}

\newtheorem{theorem}{Theorem}[section]
\newtheorem{lemma}[theorem]{Lemma}
\newtheorem{corollary}[theorem]{Corollary}

\theoremstyle{definition}

\usepackage{fancyhdr}

\newcommand{\coeff}[2]{\left[#1^{\,#2}\right]}
\newcommand{\succeqcoef}{\succeq_{\mathrm{coef}}}

\title[Hook-Length Biases in $t$-regular partitions]{Hook-Length Biases in $t$-regular partitions}

\author[M. P. Saikia]{Manjil P. Saikia\,\orcidlink{0000-0002-2997-6731}}
\address{Mathematical and Physical Sciences division, School of Arts \& Sciences, Ahmedabad University, Navrangpura, Ahmedabad 380009, Gujarat, India}
\email{manjil.saikia@ahduni.edu.in}

\author[P. Talukdar]{Prabal Talukdar\,\orcidlink{0009-0006-3713-1348}}
\address{Department of Mathematics and Statistics, Indian Institute of Science Education and Research Kolkata, 
Mohanpur, Nadia 741246, West Bengal, India
}
\email{prabaltalukdar89@gmail.com}



\begin{document}

\begin{abstract}
Recently, there has been a lot of work on combinatorial inequalities related to hook-lengths in $t$-regular partitions. In this short note, we give a proof
using generating functions for a result proved by Singh and Barman (2026) using combinatorial methods. In addition, we give an alternate proof of another result of Singh \& Barman (2024) which yields as a corollary a previously unobserved connection of hook-lengths in $t$-regular partitions with certain distinct parts partitions.
\end{abstract}

\maketitle

\section{Introduction \& Statement of Results}

A partition of a positive integer $n$ is a finite sequence of non-increasing positive integers $\lambda =(\lambda_1, \lambda_2, \ldots, \lambda_k)$ such that $\sum\limits_{i=1}^k\lambda_i=n$. The integers $\lambda_i$'s are called parts of the partition $\lambda$. There are several ways of pictorially representing partitions, one of which is via means of an \textit{Young diagram}. Such a diagram for the partition $\lambda$ is a left-justified array of boxes with the $i$-th row from the top having $\lambda_i$ boxes. For instance, the Young diagram of the partition $(4,3,1,1)$ is shown below (we ignore the numbers listed at the moment).
\begin{center}
\begin{tikzpicture}[scale=0.7]
  \tikzset{
    cell/.style={draw, minimum width=0.7cm, minimum height=0.7cm, align=center}
  }
  \node[cell] at (0,0) {7};
  \node[cell] at (1,0) {4};
  \node[cell] at (2,0) {3};
  \node[cell] at (3,0) {1};
  \node[cell] at (0,-1) {5};
  \node[cell] at (1,-1) {2};
  \node[cell] at (2,-1) {1};
  \node[cell] at (0,-2) {2};
  \node[cell] at (0,-3) {1};
\end{tikzpicture}  
\end{center}

The \textit{hook length} of a box in such a Young diagram is the sum of the number of boxes directly to the right of the box, the number of boxes directly below and $1$ for the box itself. For instance, the hook lengths of the boxes in the above Young diagram are written inside the respective boxes. In this paper, we focus on hook lengths of $t$-regular partitions. By a $t$-regular partition, we mean a partition with no part divisible by $t$. For instance, $(4,3,1,1)$ is a $5$-regular partition, while $(10,7,3,1)$ is not.

Let $t\geq 2$ be a fixed positive integer, and let $b_{t,k}(n)$ denote the number of hooks of length $k$ in all the $t$-regular partitions of $n$. Recently, there have been a lot of work related to this statistic. For instance, Ballantine et al. \cite{BallantineEtAl} studied hook lengths in $2$-regular partitions among other things, Singh and Barman \cite{SinghBarmanJNT} studied hook length biases for $2$- and $3$-regular partitions for different hook lengths. They established combinatorial inequalities such as $b_{2,2}(n)\geq b_{2,1}(n)$ for all $n>4$, etc. Very recently, Singh and Barman \cite{SinghBarmanCMB} continued their studies and proved the following results.

\begin{theorem}\cite[Theorem 1.2]{SinghBarmanCMB}\label{thm1}
    For all $n>3$, we have
\[
b_{3,2}(n)\ \ge\ b_{2,2}(n).
\]
\end{theorem}

\begin{theorem}\cite[Theorem 1.3]{SinghBarmanCMB}\label{thm2}
    For all $n\geq 0$, we have
\[
b_{3,3}(n)\ \ge\ b_{2,3}(n).
\]
\end{theorem}

\noindent Theorems \ref{thm1} and \ref{thm2} were proved combinatorially by Singh and Barman \cite{SinghBarmanCMB}. At the end of their paper, they remarked that it would be interesting to find proofs of their results using $q$-series techniques based on generating functions. We do this here, for Theorem \ref{thm1} in Section \ref{sec:proofthm1}.

In addition to Theorem \ref{thm2} above, the following results are also in need of elementary $q$-series proofs using generating functions, as their original proofs are combinatorial.

\begin{theorem}\cite[Theorem 1.4]{SinghBarmanJNT}\label{thm3}
     For all $n > 4$, we have
\[
b_{2,2}(n)\ \ge\ b_{2,1}(n).
\]
\end{theorem}

\begin{theorem}\cite[Theorem 1.2]{PJM}
     For all $n\geq 0$, we have
\[
b_{4,2}(n)\ \ge\ b_{3,2}(n).
\]
\end{theorem}

We give an alternate proof of Theorem \ref{thm3} in Section \ref{sec:thm3}, where we use combinatorial reasoning only in one step of the proof. Even though we were unable to give a purely analytic proof for Theorem \ref{thm3}, our proof yields the following previously unobserved connection with distinct parts partitions.
\begin{corollary}\label{coro}
    Let $d_3(n)$ denote the number of distinct parts partitions of $n$, where the parts are all of size greater than $2$. Then, we have for all $n>4$
    \[
    b_{2,2}(n)-b_{2,1}(n)=\sum_{i=0}^{n-5}d_3(i)-d_3(n-1),
    \]
    where we assume $d_3(0)=1$.
\end{corollary}
\noindent Corollary \ref{coro} is proved in Section \ref{sec:thm3}. It would be interesting to obtain a purely combinatorial proof of Corollary \ref{coro}.

\section{Preliminaries}\label{sec:prelim}

Throughout, we assume $|q|<1$. The $q$-Pochhammer symbol is defined as
\[
(a;q)_n:=\prod_{i=0}^{n-1}(1-aq^i) \quad \text{and} \quad (a;q)_\infty := \prod_{m\ge0}(1-aq^m).
\]
For integers $t\ge2$ and $k\ge1$, recall that $b_{t,k}(n)$ denote the number of hooks of length $k$ among all $t$-regular partitions of $n$. Set
\[
B_{t,k}(q):=\sum_{n\ge0} b_{t,k}(n)\,q^n.
\]

We have the following generating functions
\begin{equation}\label{eq:kim}
B_{2,1}(q)=\frac{1}{(q;q^2)_\infty}\Bigl(\frac{q}{1-q}-\frac{q^2}{1-q^2}\Bigr),
\end{equation}
\begin{equation}\label{eq:ballantine}
B_{2,2}(q)=\frac{1}{(q;q^2)_\infty}\Bigl(q^2+\frac{q^3}{1-q^2}+\frac{q^6}{1-q^4}\Bigr),
\end{equation}
and
\begin{equation}\label{eq:B32B22}
B_{3,2}(q)=\frac{(q^3;q^3)_\infty}{(q;q)_\infty}\Big(\frac{q^2}{1-q}+\frac{q^2}{1-q^2}-\frac{2q^3}{1-q^3}\Big).
\end{equation}
Kim \cite[Thm. 3.1]{Kim} has given an unified proof for all of these identities; \eqref{eq:ballantine} was also proved by Ballantine et al. \cite[Proposition 3.1]{BallantineEtAl}, while \eqref{eq:B32B22} was also proved by Singh and Barman \cite[Equation 1.7]{SinghBarmanJNT}.

We use the specializations 
\[
(-q;q)_\infty := \prod_{m\ge1}(1+q^m)=\frac{(q^2;q^2)_\infty}{(q;q)_\infty}
=(1+q)(1+q^2)\,(-q^3;q)_\infty\]
and\[
(-q^3;q)_\infty := \prod_{m\ge3}(1+q^m)
\]
throughout the paper without commentary. We write $A(q)\succeqcoef B(q)$ to mean $\coeff{q}{n}A(q)\ge \coeff{q}{n}B(q)$ for every $n$ (coefficient-wise order), where $[q^n]N(q)$ is the coefficient of $q^n$ in the formal power series $N(q)$. Sometimes we say the coefficients of $N(q)$ are nonnegative for all $n\geq k$, when we mean that the coefficients of $q^n$ in the power series $N(q)$ are nonnegative when $n\geq k$.

We need two classical identities, the first due to Euler \cite[p. 5]{gea1}:
\begin{equation}\label{euler}
    (-q;q)_\infty=\frac{1}{(q;q^2)_\infty},
\end{equation}
and the second due to Sylvester \cite[p. 281]{Sylvester}:
\begin{align}\label{sylvester}
(-xq;q)_{\infty} = \sum_{n\geq0}\frac{(-xq;q)_{n}}{(q;q)_n}(1+xq^{2n+1})x^nq^{n(3n+1)/2}. 
\end{align}
Both of these identities have classical proofs using $q$-series. Finally, We need the two lemmas below.
\begin{lemma}\label{lem:coef-positivity}
Let
\[
H(q)=\sum_{n\ge 0} h_n q^n,
\]
where $(h_n)_{n\ge0}$ is eventually nondecreasing; that is, there exists $N_0$ such that
$h_{n+1}\ge h_n$ for all $n\ge N_0$.
Let $r\in \mathbb{N}$  and  $J\subseteq \{0\}\cup\mathbb{N}$ be a set with $r\notin J$. Also, let
\[
E(q)=\Big(\sum_{j\in J} \alpha_j q^j\Big)-c\,q^{r},
\quad\text{with } \alpha_j\ge0 \text{ for } j\in J \text{ and } c\ge0.
\]
If there exists $s\in J$ with $s<r$  with $\alpha_s\geq c$, then for all $n\ge N_0+r$,
\[
[q^n]\big(E(q)H(q)\big)\ge 0 .
\]
\end{lemma}

\begin{proof}
For $n\ge N_0+r$ we have
\[
[q^n]\big(E(q)H(q)\big)
= \sum_{j\in J} \alpha_j\, h_{n-j} - c\, h_{n-r}
\ge \alpha_s h_{n-s} - c\, h_{n-r}
\ge c\, h_{n-s} - c\, h_{n-r}
= c\big(h_{n-s}-h_{n-r}\big).
\]
Since $s<r$, we have $n-s>n-r$. For $n\ge N_0+r$ both $n-s$ and $n-r$ are at least $N_0$, 
and because $(h_n)$ is eventually nondecreasing, $h_{n-s}\ge h_{n-r}$. Hence
$c\big(h_{n-s}-h_{n-r}\big)\ge0$, proving the claim.
\end{proof}

\begin{lemma}\label{lem:d3-monotone}
Let us write
\[
(-q^{3};q)_\infty \;=\; \prod_{m\ge3}(1+q^m) \;:=\; \sum_{n\ge0} d_3(n)\,q^n,
\]
where clearly $d_3(n)\ge 0$ for all $n$. Combinatorially, $d_3(n)$ is the number of
partitions of $n$ into distinct parts, all of whose sizes are at least $3$.
Then, for all $n\ge2$,
\[
d_3(n)\;\ge\; d_3(n-1).
\]
Here we assume $d_3(0)=1$.
\end{lemma}

\noindent We can prove Lemma \ref{lem:d3-monotone} quite easily using combinatorial reasoning, but we give an analytic proof below because we wish to give a completely analytic proof of Theorem \ref{thm1}, and Lemma \ref{lem:d3-monotone} is used in the proof.

\begin{proof}[Proof of Lemma \ref{lem:d3-monotone}]
    Notice that
    \begin{equation}\label{eq:lem-1}
        (1-q)(-q^3;q)_\infty=\sum_{n\geq 1}(d_3(n)-d_3(n-1))q^n+d_3(0).
    \end{equation}
    To prove the lemma, it is enough to prove that the left-hand side of equation \eqref{eq:lem-1} has nonnegative coefficients for all $n\geq 2$. Let us rewrite, the left-hand side as follows (where we use \eqref{euler} in the first step)
    \begin{align}\label{eq:syl}
        (1-q)(-q^3;q)_\infty&=\frac{(1-q)(1-q^2)}{(1+q)(1-q^4)}\frac{1}{(q;q^2)_\infty}\nonumber\\
        &=\frac{1-q}{1-q^4}\frac{1}{(q^3;q^2)_\infty}\nonumber\\
        &=\left(\frac{1}{(q^3;q^2)_\infty}\frac{q^2(1-q)}{(1-q^2)^2}\right)\frac{(1-q^2)}{q^2(1+q^2)}.
    \end{align}

We put $x = 1$ in equation \eqref{sylvester} to obtain the following
\begin{align*}
 q^3+q^5+\frac{1}{(q^3;q^2)_{\infty}}\frac{q^2(1-q)}{(1-q^2)^2}
& = q^3(1+q^2)+\frac{q^2(1-q)^2}{(1-q^2)^2}(-q;q)_{\infty}\\
& = q^3(1+q^2)+\frac{q^2}{(1+q)^2}\sum_{n\geq0}\frac{(-q;q)_{n}}{(q;q)_{n}}(1+q^{2n+1})q^{\frac{3n^2+n}{2}}\\
& = q^3(1+q^2)+\frac{q^2}{(1+q)}+\frac{q^4(1+q^3)}{(1-q^2)}\\ &\quad +\frac{q^2}{(1+q)^2}\sum_{n\geq2}\frac{(-q;q)_{n}}{(q;q)_{n}}(1+q^{2n+1})q^{\frac{3n^2+n}{2}}\\
& = \frac{q^2(1+q^2)}{(1-q^2)}+\frac{q^2}{(1+q)^2}\sum_{n\geq2}\frac{(-q;q)_{n}}{(q;q)_{n}}(1+q^{2n+1})q^{\frac{3n^2+n}{2}},
\end{align*}
which gives
\begin{align}\label{eq:syl-2}
 \frac{1}{(q^3;q^2)_{\infty}}\frac{q^2(1-q)}{(1-q^2)^2} = -q^3-q^5+\frac{q^2(1+q^2)}{(1-q^2)}+\frac{q^2}{1-q^2}\sum_{n\geq2}\frac{(-q^2;q)_{n-1}}{(q^2;q)_{n-1}}(1+q^{2n+1})q^{\frac{3n^2+n}{2}}.
\end{align}

Using \eqref{eq:syl-2} in \eqref{eq:syl}, we obtain after some simplification
\begin{align}\label{eq:syl-3}
    (1-q)(-q^3;q)_\infty&=1-q+q^3+\sum_{n\geq 2}\frac{(-q^3;q)_{n-2}}{(q^2;q)_{n-1}}(1+q^{2n+1})q^{(3n^2+n)/2}.
\end{align}
Notice that all the coefficients of $q^n$ in the sum on the right-hand side of \eqref{eq:syl-3} are positive, so combining \eqref{eq:syl} -- \eqref{eq:syl-3} we have proved the result for all $n\geq 3$. The result clearly holds for $n=2$, as we know that $d_3(2)=0=d_3(1)$.
\end{proof}

\section{Proof of Theorem \ref{thm1}}\label{sec:proofthm1}

Let
\[
G(q) := B_{3,2}(q) - B_{2,2}(q) 
= \sum_{n \geq 0} \big( b_{3,2}(n) - b_{2,2}(n) \big) q^n
= \sum_{n \geq 0} G_n q^n.
\]
Thus, to prove Theorem \ref{thm1} it suffices to show that $G_n \geq 0$ for all $n \geq 4$.

Simplifying \eqref{eq:ballantine}, we get
\[
B_{2,2}(q)=\sum_{n\ge0} b_{2,2}(n)q^n
=(-q^3;q)_\infty (1+q)(1+q^2)\Bigl(q^2+\frac{q^3}{1-q^2}+\frac{q^6}{1-q^4}\Bigr).
\]
A short computation shows that
\begin{equation}
\ B_{2,2}(q)=(-q^3;q)_\infty\frac{q^2+2q^3+q^4+q^5+q^6}{1-q^2}\ .
\end{equation}

We denote
\begin{multline}\label{eq:CRP}
C(q):=\frac{q^2}{1-q}+\frac{q^2}{1-q^2}-\frac{2q^3}{1-q^3},\quad
R(q):=\frac{q^2+2q^3+q^4+q^5+q^6}{1-q^2}, \quad \text{and}\\
P(q):=\frac{(q^3;q^3)_\infty}{(q;q)_\infty}=\prod_{m\ge1}(1+q^m+q^{2m}).
\end{multline}
Note that $P(q)$ has nonnegative coefficients. We set
\[
F(q):=(1-q^6)\,G(q).
\]
We will prove that $F(q)$ has nonnegative coefficients $[q^n]$ for all $n\geq 4$, and connect the coefficients of $F(q)$ to the coefficients of $G(q)$.

For ease of notation, we rewrite $F(q)$ using \eqref{eq:ballantine} and \eqref{eq:B32B22} as follows
\begin{align*}
    F(q)&=(1-q^6)\left[\frac{(q^3;q^3)_\infty}{(q;q)_\infty}\Big(\frac{q^2}{1-q}+\frac{q^2}{1-q^2}-\frac{2q^3}{1-q^3}\Big)-\frac{1}{(q;q^2)_\infty}\Bigl(q^2+\frac{q^3}{1-q^2}+\frac{q^6}{1-q^4}\Bigr)\right]\nonumber\\
    &=P(q)\,P_C(q)-(-q^3;q)_\infty\,P_R(q),\label{eq:Fdef}
\end{align*}
where $C(q), R(q)$ and $P(q)$ are given by \eqref{eq:CRP}, and
\[
P_C(q):=(1-q^6)C(q),\qquad P_R(q):=(1-q^6)R(q).
\]
Simplifying this further, we have
\begin{equation}\label{eq:PR}
 P_R(q)=q^2+2q^3+2q^4+3q^5+3q^6+3q^7+2q^8+q^9+q^{10}.
\end{equation}

Define the polynomial
\[
M(q):=(1+q+q^2)(1+q^2+q^4).
\]

\begin{lemma}\label{lem:Phi}
Set
\[
\Phi(q):=M(q)\,P_C(q)-P_R(q)-q^3.
\]
Then $\coeff{q}{n}\Phi(q)\ge0$ for all $n\ge4$.
\end{lemma}
\begin{proof}
This is easily checked by expanding:
\[
\Phi(q)=q^2-2q^3+3q^4+5q^6+2q^7+5q^8+4q^9+3q^{10}+3q^{11}+q^{12}+q^{13}.
\]
\end{proof}

\begin{lemma}\label{eq:Pge}
    We have,
    \[ P(q)\ \succeqcoef\ M(q)\cdot \prod_{m\ge3}(1+q^m)=M(q)\cdot (-q^3;q)_\infty .\]
\end{lemma}

\begin{proof}
We split $P(q)$ as
\[
P(q)=(1+q+q^2)(1+q^2+q^4)\prod_{m\ge3}(1+q^m+q^{2m})=M(q)\cdot \prod_{m\ge3}(1+q^m+q^{2m}).
\]
Notice that the expansion of $\prod_{m\geq 3}(1+q^m+q^{2m})$ will contain a copy of $\prod_{m\geq 3}(1+q^m)$, which gives us the result.
\end{proof}

We now define
\[
S(q):=P(q)-M(q)(-q^3;q)_\infty.
\]
Then, 
    we have
    \[
    [q^n]S(q)= 0, \quad \text{for all $0\leq n\leq 5$},
    \]
    since
    \begin{align*}
    S(q)&=\prod_{m\geq 1}(1+q^m+q^{2m})-(1+q+q^2)(1+q^2+q^4)\prod_{m\geq 3}(1+q^m),
\end{align*}
and by a short computation, it is clear that $[q^n]S(q)=0$ for $0\leq n \leq 5$. In addition, 
     \[
    [q^n]S(q) \ge  0, \quad \text{for all $ n\geq 6$},
    \]   
which is just a restatement of Lemma \ref{eq:Pge}. Also note
\[
F(q)=P(q)P_C(q)-(-q^3;q)_\infty P_R(q).
\]
\begin{lemma}\label{lem:fn}
    We have,
    \[
    [q^n]F(q)\ge0\ \text{ for all }n\ge4.
    \]
\end{lemma}

\begin{proof}
    From the definition of $F(q)$, we have
\begin{align*}
    F(q)&=(M(q)(-q^3;q)_\infty+S(q))P_C(q)-(-q^3;q)_\infty P_R(q)\\
    &=(-q^3;q)_\infty (M(q)P_C(q)-P_R(q))+S(q)P_C(q).
\end{align*}

We first claim that $[q^n](-q^3;q)_\infty (M(q)P_C(q)-P_R(q))\geq 0$ for $n\geq 4$. Notice,
\begin{align*}
(-q^3;q)_\infty\big(M(q)P_C(q)-P_R(q)\big)
&=\prod_{m\ge 3}(1+q^m)\big(\Phi(q)+q^3\big)\\[4pt]
&=\prod_{m\ge 3}(1+q^m)\big(q^2+\Phi_{\ge4}(q)-q^3\big)\\[-2pt]
&\quad \text{(where }\Phi_{\ge4}(q)\text{ is the polynomial with}\\& \quad \text{ all the degree }\ge 4\text{ terms of $\Phi(q)$)}\\[6pt]
&\ge \prod_{m\ge 3}(1+q^m)\big(q^2+3q^4+q^6-q^3\big).
\end{align*}
So, to prove our claim it suffices to show that $\prod_{m\geq 3}(1+q^m) (q^2+3q^4+q^6-q^3)$ has nonnegative coefficients from the degree $4$ term onward.

Let
\[
H(q)=(-q^{3};q)_\infty
=\prod_{m\ge3}(1+q^m)
=\sum_{n\ge0} d_3(n)\,q^n,
\qquad
E(q)= q^{2}+3q^{4}+q^{6}-q^{3}.
\]
By Lemma~\ref{lem:d3-monotone} we know $d_3(n+1)\ge d_3(n)$ for all $n\ge1$, so the coefficient
sequence of $H(q)$ is eventually nondecreasing with \(N_0=1\).
Applying Lemma~\ref{lem:coef-positivity} with \(r=3\), \(s=2\) and \(c=1\) 
\[
[q^n]\big(E(q)H(q)\big)\;\ge\;0
\qquad\text{for all } n\ge N_0+r = 4.
\]
So, we have proved our claim, that 
\begin{equation}\label{eq:f1}
    [q^n](-q^3;q)_\infty(M(q)P_C(q)-P_R(q))\geq 0, \quad n\geq 4.
\end{equation}

All that remains to show now is that 
\begin{equation}\label{eq:f3}
[q^n]S(q)P_C(q)\geq 0, \quad n\geq 0. 
\end{equation}
It is easy to see computationally that this is true for $8\geq n \geq 0$. We claim that this holds for $n\geq 9$. We have,
\[
S(q)P_C(q)=M(q)\left(\prod_{m\geq 3}(1+q^m+q^{2m})-\prod_{m\geq 3}(1+q^m)\right)P_C(q).
\]
By a direct computation, we have
\[
M(q)P_C(q)=\left(q^3+q^2+q+2\right) \left(q^5+q^3+q\right)^2,
\]which clearly has nonnegative coefficients. Additionally, the expansion of $\prod_{m\geq 3}(1+q^m+q^{2m})$ will contain a copy of $\prod_{m\geq 3}(1+q^m)$ and hence the coefficients of $\prod_{m\geq 3}(1+q^m+q^{2m})-\prod_{m\geq 3}(1+q^m)$ are nonnegative. Hence, we have
\begin{equation*}
    [q^n]S(q)P_C(q)\geq 0, \quad n\geq 9.
\end{equation*}

Combining \eqref{eq:f1} and \eqref{eq:f3} we have the required result.
\end{proof}

We are now ready to complete the proof of Theorem \ref{thm1}. Recall that
\[
G(q)=B_{3,2}(q)-B_{2,2}(q)=\sum_{n\ge0}G_n q^n,
\]
and
\[
F(q)=(1-q^6)\,G(q):=\sum_{n\ge0}F_n q^n.
\]
We have
\[
F(q)=(1-q^6)G(q)=G(q)-q^6G(q).
\]
Now
\[
q^6G(q)=\sum_{n\ge0}G_n q^{n+6}
      =\sum_{m\ge6}G_{m-6}q^{m}
      =\sum_{n\ge0}G_{n-6}q^{n},
\]
where we interpret \(G_{n-6}=0\) for \(n<6\). Hence
\[
F(q)=\sum_{n\ge0}\big(G_n-G_{n-6}\big)q^n.
\]
Thus, for all integers \(n\geq 0\),
\begin{equation}\label{eq:rec6}
G_n-G_{n-6}=F_n.
\end{equation}

In Lemma \ref{lem:fn} we proved
\begin{equation}\label{eq:Fge0}
F_n\ge0 \qquad (n\ge4).
\end{equation}
From \eqref{eq:rec6} and \eqref{eq:Fge0} it follows immediately that
\begin{equation}\label{eq:monotone}
G_n\ge G_{n-6}\qquad(n\ge4).
\end{equation}
We compute $G_n$ for $n=4,5,6,7,8,9$ directly from the definitions of $B_{3,2}$ and $B_{2,2}$. This yields 
\begin{equation}\label{eq:baseblock}
(G_4,G_5,G_6,G_7,G_8,G_9)=(3,1,5,5,11,13).
\end{equation}

Hence, for each $r\in\{4,5,6,7,8,9\}$, the sequence $(G_{r+6k})_{k\ge 0}$ is nondecreasing by \eqref{eq:monotone}. Since \eqref{eq:baseblock} holds, it follows that $G_n\ge 0$ for all
$n\ge 4$. This completes the proof of Theorem \ref{thm1}.

\section{Proofs of Theorem \ref{thm3} and Corollary \ref{coro}}\label{sec:thm3}

\begin{proof}[Proof of Theorem \ref{thm3}]
From \eqref{eq:kim} and \eqref{eq:ballantine}, we have
\begin{align}
    B_{2,2}(q)-B_{2,1}(q)&=(-q;q)_\infty\left( \frac{q^2+q^3+q^5}{1-q^4}-\frac{q}{1-q^2}\right)\nonumber\\ 
    &=(-q;q)_\infty \left(\frac{q^2+q^5-q}{1-q^4}\right)\nonumber\\ 
    &=(-q^3;q)_\infty\left(\frac{q^2+q^5-q}{1-q}\right)\nonumber\\ 
    &=(-q^3;q)_\infty\left(q^2\sum_{n\geq 0}q^n+q^5\sum_{n\geq 0}q^n-q\sum_{n\geq 0}q^n\right)\nonumber\\
    &=(-q^3;q)_\infty\left(-q+\sum_{n\geq 5}q^n\right)\label{b22}.
\end{align}

Recall
    \[
    (-q^3;q)_\infty=\sum_{n\geq 0}d_3(n)q^n,
    \]
    where clearly $d_3(n)\geq 0$ for all $n$ (see the statement of Lemma \ref{lem:d3-monotone}). We now rewrite \eqref{b22} as follows
    \begin{align}
      B_{2,2}(q)-B_{2,1}(q)&=\left(\sum_{n\geq 0}d_3(n)q^n\right)\left(\sum_{n\geq 5}q^n\right) -\sum_{n\geq 0}d_3(n)q^{n+1} \nonumber\\
      &=\left(\sum_{n\geq 0}d_3(n)q^n\right)\left(\sum_{n\geq 0}q^n\right)-\left(\sum_{n\geq 0}d_3(n)q^n\right)(1+q+q^2+q^3+q^4)\nonumber\\
      &\quad -\sum_{n\geq 0}d_3(n)q^{n+1}\nonumber\\
      &=\sum_{n\geq 0}c(n)q^n-\sum_{n\geq 0}d_3(n)q^{n}-2\sum_{n\geq 0}d_3(n)q^{n+1}-\sum_{n\geq 0}d_3(n)q^{n+2}\nonumber\\
      &\quad -\sum_{n\geq 0}d_3(n)q^{n+3}-\sum_{n\geq 0}d_3(n)q^{n+4},\label{b2-2}
    \end{align}
    where
    \[
    c(n)=\sum_{i=0}^nd_3(i).
    \]
Simplifying \eqref{b2-2}, we obtain
\begin{align}
     B_{2,2}(q)-B_{2,1}(q)&=\sum_{n\geq 0}(c(n)-d_3(n))q^n-2\sum_{n\geq 1}d_3(n-1)q^{n}-\sum_{n\geq 2}d_3(n-2)q^{n}\nonumber\\
     &\quad -\sum_{n\geq 3}d_3(n-3)q^{n}-\sum_{n\geq 4}d_3(n-4)q^{n}\nonumber\\
     &=\sum_{n\geq 4}\left(c(n)-d_3(n)-2d_3(n-1)-d_3(n-2)-d_3(n-3)-d_3(n-4)\right)q^n\nonumber\\
     &\quad +\sum_{n=0}^3(c(n)-d_3(n))q^n-2q-q^2-q^3,\label{b22-3}
\end{align}
where we have used the fact that $d_3(0)=1,d_3(1)=d_3(2)=0$.

From \eqref{b22-3}, it is clear that to prove Theorem \ref{thm3} we need to show the following
\begin{equation}\label{b22-4}
c(n)-d_3(n)-2d_3(n-1)-d_3(n-2)-d_3(n-3)-d_3(n-4) \geq 0 \quad \text{for all}~n>4.
\end{equation}
Simplifying \eqref{b22-4}, using the value of $c(n)$, this is equivalent to showing the following
\begin{equation}\label{b22-5}
    \sum_{i=0}^{n-5}d_3(i)-d_3(n-1)\geq 0\quad \text{for all}~n>4.
\end{equation}
We now prove \eqref{b22-5} using a combinatorial argument. Let $D_3(n)$ denote the set of all partitions of $n$ with distinct parts greater than $2$. Thus, $|D_3(n)|=d_3(n)$. We notice that among all partitions in $D_3(n-1)$ there are no partitions with any part size equal to $n-2$ or $n-3$ due to the restriction on parts. We map the unique partition $(n-1)\in D_3(n-1)$ to the unique partition in $\phi \in D_3(0)$. For every other partition $\lambda=(\lambda_1, \lambda_2, \ldots, \lambda_k)\in D_3(n-1)$, we just delete the largest part $\lambda_1$ and map it to a partition in one of the $D_3(i)$'s, where $3\leq i\leq n-5$. This gives us an injection, and we are done.
\end{proof}

\begin{proof}[Proof of Corollary \ref{coro}]
    The result follows from equations \eqref{b22-3}, \eqref{b22-4} and \eqref{b22-5}.
\end{proof}

\section*{Acknowledgments}

The authors thank the anonymous referee for several helpful comments which improved the exposition. The first author is partially supported by a Start-Up Grant from Ahmedabad University (Reference No. URBSASI24A5).

\newcommand{\etalchar}[1]{$^{#1}$}

\end{document}